# HIGH ORDER TENSOR-TRAIN-BASED SCHEMES FOR HIGH-DIMENSIONAL MEAN FIELD GAMES*

ELISABETTA CARLINI† AND LUCA SALUZZI‡

**Abstract.** We introduce a fully discrete scheme to solve a class of high-dimensional Mean Field Games systems. Our approach couples semi-Lagrangian (SL) time discretizations with Tensor-Train (TT) decompositions to tame the curse of dimensionality. By reformulating the classical Hamilton–Jacobi–Bellman and Fokker–Planck equations as a sequence of advection–diffusion–reaction subproblems within a smoothed policy iteration, we construct both first and second order in time SL schemes. The TT format and appropriate quadrature rules reduce storage and computational cost from exponential to polynomial in the dimension. Numerical experiments demonstrate that our TT-accelerated SL methods achieve their theoretical convergence rates, exhibit modest growth in memory usage and runtime with dimension, and significantly outperform grid-based SL in accuracy per CPU-second.

## 1. Introduction.

Mean Field Games (MFGs) theory provides a powerful mathematical framework for modeling the collective behavior of a large number of rational agents who interact through their empirical distribution [33, 36]. The resulting systems couple a Hamilton-Jacobi-Bellman (HJB) equation, describing the optimal control of a representative agent, with a Fokker-Planck (FP) equation, governing the evolution of the population density. From a numerical point of view, solving such coupled systems remains highly challenging, especially in high-dimensional settings where traditional grid-based discretizations suffer from the curse of dimensionality.
Many works have focused on developing numerical methods for Mean Field Games, ranging from finite difference, semi-Lagrangian schemes to finite element methods [2, 1, 18, 45]. More recently, several contributions have addressed the solution of high-dimensional Mean Field Games, with most approaches relying on neural networks. In particular, neural network methods have been developed for advection-diffusion problems and high-dimensional optimal control [28, 44]. Notable examples include deep reinforcement learning and physics-informed neural networks for MFGs [52, 19, 29, 4, 5], which provide flexible and scalable approximations in very high-dimensional settings. Another important aspect concerns the development of second order schemes, which play a crucial role in improving accuracy without excessive computational cost. Recent developments on high order numerical methods for Mean Field Games include finite difference $\theta$-scheme, multiscale approach and SL-Lagrange Galerkin type schemes [50, 6, 38, 7]. However, high order numerical schemes that remain efficient in high-dimensional settings are still largely unexplored.

A third aspect of the approximation is related to the development of efficient algorithms for solving the nonlinear MFGs system, either in its continuous or discretized

*

**Funding:** This work was funded by Italian Ministry of Instruction, University and Research (MIUR) (PRIN Project 2022 2022238YY5, “Optimal control problems: analysis, approximation”) and INdAM-research group GNCS (CUP E53C24001950001, “Metodi avanzati per problemi di Mean Field Games ed applicazioni”) form. Several procedures have been recently proposed, including policy iteration and Newton-type algorithms [3, 9, 10, 54]. Building on the policy iteration algorithm introduced in [54], we propose a second order semi-Lagrangian (SL) discretization combined with Tensor-Train (TT) techniques. An increasing interest has emerged in the development of low-rank and tensor-based techniques to mitigate the complexity of high-dimensional PDEs arising in control and game theory. Among these approaches, the Tensor-Train decomposition [46] has emerged as an effective tool for representing high-dimensional functions and operators in a compact, structured format. Tensor-based discretizations have been successfully employed in the context of high-dimensional FP equations [25, 20, 22] and HJB equations [23, 49, 24]. In particular, for FP equations, the solution often exhibits smoothness that can be effectively exploited by low-rank tensor approximations, leading to highly efficient and accurate numerical schemes.

†Department of Mathematics, Sapienza, University of Rome, Italy (carlini@mat.uniroma1.it)

‡Department of Mathematics, Universitá di Ferrara, Italy (luca.saluzzi@unife.it)

In this work, we move toward the development of accurate and computationally feasible high order methods for a class of high-dimensional MFGs problem, particularly addressing second order problems and first order problems with smooth solutions. We combine the stability properties of SL discretizations with the efficiency of TT representations. SL schemes allow us to derive an unconditionally stable semi-discrete in time approximations of the linearized HJB and FP equations, inherent in the Policy Iteration algorithm for MFGs models. TT reconstruction allows us to build a fully discrete scheme on a tensor grid with low computational cost. By combining these two techniques, we achieve a scalable and memory-efficient numerical scheme that mitigates the exponential growth of computational complexity. We propose a new second order semi-Lagrangian scheme based on quadrature nodes, whose computational cost grows only polynomially with the dimension. A potential drawback of the proposed method is that, in order to achieve polynomial cost, negative quadrature weights must be employed. Nevertheless, we perform extensive numerical simulations to demonstrate the robustness and accuracy of the scheme, and we show that it is asymptotically positivity-preserving. We discuss some theoretical properties of the method, such as numerical accuracy, and demonstrate its effectiveness through computational experiments on high-dimensional test cases.

The paper is organized as follows. In Section 2, we recall the policy iteration framework for Mean Field Games and introduce the proposed second order semi-Lagrangian discretization in time. In Section 3, we discuss the approximation of the expectation term arising in the semi-Lagrangian formulation and introduce the quadrature strategy that ensures second order accuracy while maintaining polynomial computational complexity. In Section 4, we provide a brief overview of TT decomposition and its use for efficient representation of high-dimensional value and density functions. In Section 5, we combine these ingredients and present the Tensor-Train-based semi-Lagrangian method, detailing its algorithmic structure. Finally, in Section 6, we report some numerical experiments that validate the accuracy, robustness, and scalability of the proposed approach, and compare its performance with standard SL schemes.

**2. Semi-Lagrangian scheme based on Policy Iteration Method.** We consider the following second order Mean Field Games system

$$
\begin{cases} -\partial_t u - \nu\Delta u + H(\nabla u) = F(x, m(t)) & \text{in } \mathbb{R}^d \times [0,T], \\ u(x,T) = G(x, m(T)) & \text{in } \mathbb{R}^d \ , \\ \partial_t m - \nu\Delta m - \operatorname{div}(m\nabla_p H(\nabla u)) = 0, & \text{in } \mathbb{R}^d \times [0,T] \\ m(x,0) = m_0(x), & \text{in } \mathbb{R}^d \ , \end{cases} \tag{2.1}
$$

where $T > 0$, and $\nu \geq 0$ is a constant diffusion parameter.

Let us remark that if $\nu > 0$ and under suitable assumptions, namely, the convexity of $H : \mathbb{R}^d \to \mathbb{R}$, the smoothness and monotonicity of $F$ and $G$ and the fact that $m_0$ is a probability density function, the system (2.1) admits a unique classical solution. On the other hand, when $\nu = 0$, the system generally fails to admit a classical solution, and well-posedness can be proved in the class of viscosity/distributional solutions (see, e.g., [13, 35, 36]). The notation $F(x, m(t))$ is used to denote a (possibly nonlocal) coupling term depending on space and on the population distribution $m(t)$. In particular, this emphasizes that the coupling may depend on the entire density $m(t)$, for instance through integral operators or global quantities such as moments of the distribution, as illustrated in Section 6.3. In this framework, the coupling functions are defined as $F : \mathbb{R}^d \times \mathcal{P}_1(\mathbb{R}^d) \to \mathbb{R}$ and $G : \mathbb{R}^d \times \mathcal{P}_1(\mathbb{R}^d) \to \mathbb{R}$, where $\mathcal{P}_1(\mathbb{R}^d)$ denotes the metric space of Borel probability measures on $\mathbb{R}^d$ with finite first order moment. Conversely, in the case of a local coupling, the dependence on the distribution is pointwise in space. Here, the coupling functions are defined as $F : \mathbb{R}^d \times \mathbb{R}^+ \to \mathbb{R}$ and $G : \mathbb{R}^d \times \mathbb{R}^+ \to \mathbb{R}$, and act directly on the value of the density at each spatial point. Accordingly, the notation $F(x, m(t))$ is understood as $F(x, m(x,t))$, as in the example of Section 6.2, where $F(x, m(x,t)) = \gamma \ln(m(x,t))$. Nonlocal couplings are particularly relevant in models where agents interact through aggregate quantities, such as average prices, first moments of the distribution, or interest rates, rather than through the local population density. From an analytical perspective, nonlocal couplings often act as smoothing operators, which facilitates the existence of classical solutions. In contrast, local couplings, especially in the deterministic regime $\nu = 0$, typically require a more delicate analysis involving weak solution frameworks; see, for example, [15].

Let us define the Lagrangian as the Legendre transform of $H$:

$$
L(q) = \sup_{p \in \mathbb{R}^d} \{p \cdot q - H(p)\}.
$$

Since $\|\nabla u\|_{L^\infty(\mathbb{R}^d \times [0,T])} \leq R$, for some $R > 0$ (see, e.g., [11, 21, 18]), the value function $u$ also satisfies

$$
-\partial_t u(x,t) - \nu\Delta u(x,t) + \sup_{|q| \leq R} \{q \cdot \nabla u(x,t) - L(q)\} = F(x, m(t)). \tag{2.2}
$$

To solve system (2.1), two main algorithms have been proposed in the literature: the Fictitious Play method [14] and the Smoothed Policy Iteration algorithm [54]. The Fictitious Play method is based on agents repeatedly computing a best response to the empirical average of past population distributions. At each iteration, it requires the solution of a full optimal control problem, consisting of a backward HJB equation followed by a forward FP equation. In contrast, Smoothed Policy Iteration involves the solution of a forward FP equation together with a backward HJB equation that is linearized around the current smoothed policy. As a consequence, the cost of the policy evaluation step is significantly reduced compared to the fully nonlinear HJB solve required by Fictitious Play. Both algorithms are known to converge globally for second-order potential Mean Field Games, starting from arbitrary initial guesses. For first-order Mean Field Games, however, due to the lack of regularity of solutions, convergence results are currently available only for the Fictitious Play algorithm. Smoothed Policy Iteration offers a more efficient alternative specifically for second-order Mean Field Games. Since we are interested here in second order problems and first order problems with smooth solutions, we consider a Smoothed Policy Iteration method proposed in [54] to solve (2.1) with the HJB equation replaced by (2.2). For completeness, we provide a self-contained description of the method in the following.

**Smoothed Policy Iteration algorithm (SPI):** Given $R > 0$, $(\delta_n)_{n\in\mathbb{N}} \in (0,1)$, and a measurable vector field $q^{(0)} : \mathbb{R}^d \times [0,T] \to \mathbb{R}^d$ with $\|q^{(0)}\|_{L^\infty(\mathbb{R}^d\times[0,T];\mathbb{R}^d)} \le R$ and $\|\mathrm{div}(q^{(0)})\|_{L^\infty(\mathbb{R}^d\times[0,T])} \le R$. Let $\overline{q}^{(0)} = q^{(0)}$ and iterate for $n \ge 0$:

**(i) Generation of the distribution from the current policy**

$$\begin{cases} \partial_t m^{(n)} - \nu\Delta m^{(n)} - \overline{q}^{(n)} \cdot \nabla m^{(n)} - m^{(n)}\mathrm{div}(\overline{q}^{(n)}) = 0 & \text{in } \mathbb{R}^d \times [0,T], \\ m^{(n)}(x,0) = m_0(x) & \text{in } \mathbb{R}^d. \end{cases}$$

**(ii) Policy evaluation**

$$\begin{cases} -\partial_t u^{(n)} - \nu\Delta u^{(n)} + \overline{q}^{(n)} \cdot \nabla u^{(n)} - L(\overline{q}^{(n)}) = F(x, m^{(n)}) & \text{in } \mathbb{R}^d \times [0,T], \\ u^{(n)}(x,T) = G(x, m^{(n)}(x,T)) & \text{in } \mathbb{R}^d. \end{cases}$$

**(iii) Policy update**

$$q^{(n+1)}(x,t) = \arg\max_{|q|\le R} \left\{ q \cdot \nabla u^{(n)}(x,t) - L(q) \right\} \qquad \text{in } \mathbb{R}^d \times [0,T].$$

**(iv) Smoothing**

$$\overline{q}^{(n+1)} = (1-\delta_n)\overline{q}^{(n)} + \delta_n q^{(n+1)}.$$

In the non-local case, under additional regularity assumptions on the data, particularly when the costs derive from a potential, [54, Theorem 3.5] shows that the SPI algorithm with $\delta_n = \frac{2}{n+2}$ converges to the solution of (2.1). Furthermore, under the Lasry–Lions monotonicity condition on the running cost F and terminal cost $G$, namely, for any $m, m' \in \mathcal{P}_1(\mathbb{R}^d)$ and $z = F, G$,

$$\int_{\mathbb{R}^d} \left(z(x,m) - z(x,m')\right) d(m-m')(x) \ge 0,$$

the convergence of the algorithm is global.

In what follows, we assume that the Hamiltonian $H$ is of class $C^2$ and strictly convex, that $u^{(n)}$ is a classical solution, and that the optimal response remains in the interior of the admissible set. Under these assumptions, see [37], the policy update admits the explicit form and step **(iii)** of the algorithm reads:

**(iiib)**

$$q^{(n+1)}(x,t) = \nabla_p H\big(\nabla u^{(n)}(x,t)\big).$$

**2.1. Second order Semi-Lagrangian approximation in time.** Solving the Mean Field Games system numerically via the Policy Iteration method leads to an iterative procedure consisting of two steps: (i) the approximation of an advection–diffusion equation with a source term $f$, and (ii) the approximation of an advection–diffusion–reaction equation. Consequently, this section is devoted to the approximation in time of these two types of equations through the use of the Feynman-Kac formula and SL schemes.

Let $N_T \in \mathbb{N}^+$ be the number of time steps, let $\Delta t = T/N_T$ be the time step, and $t_k = k\Delta t$ for all $k = 0, \dots, N_T$.

**Feynman–Kac formula for forward FP equation.** Suppose $r : \mathbb{R}^d \times [0,T] \to \mathbb{R}$, $q : \mathbb{R}^d \times [0,T] \to \mathbb{R}^d$ are given functions, and let us write a general forward linear advection–diffusion-reaction equation, as appears in steps **(i)**, in the form:

$$\partial_t m - \nu \Delta m - q \cdot \nabla m - r\, m = 0. \tag{2.3}$$

The Feynman–Kac representation provides a probabilistic expression for the solution of (2.3) over a time interval $[t_k, t_{k+1}]$. It can be written as

$$m(x,t_{k+1}) = \mathbb{E}\left[ m(X^-_{x,t_{k+1}}(t_k), t_k)\, \exp\left( \int_{t_k}^{t_{k+1}} r(X^-_{x,t_{k+1}}(\tau), \tau)\, d\tau \right) \middle| X^-_{x,t_{k+1}}(t_{k+1}) = x \right], \tag{2.4}$$

where $X^-_{x,t_{k+1}}(\cdot)$ denotes a stochastic process initiated at position $x$ at time $t_{k+1}$ evolving backward in time. The process is governed by the drift term $q$ and by a Wiener process $\sqrt{2\nu} dW$. For standard references on the Feynman–Kac formula, see [27, 43]. This formula can be discretized in time in order to heuristically derive a semi-discrete semi-Lagrangian-type scheme for (2.3), as in [26, 34]. Applying a trapezoidal rule in $[t_k, t_{k+1}]$, we get

$$\int_{t_k}^{t_{k+1}} r(X^-_{x,t_{k+1}}(\tau), \tau)\, d\tau \simeq \frac{\Delta t}{2}\left( r(X^-_{x,t_{k+1}}(t_k), t_k) + r(x, t_{k+1}) \right),$$

which allows us to approximate (2.4) as

$$m(x,t_{k+1}) \simeq \mathbb{E}\left[ m\big(X^-_{x,t_{k+1}}(t_k), t_k\big)\, \exp\left( \tfrac{\Delta t}{2} \left( r\big(X^-_{x,t_{k+1}}(t_k), t_k\big) + r\big(x, t_{k+1}\big) \right) \right) \middle| X^-_{x,t_{k+1}}(t_{k+1}) = x \right]. \tag{2.5}$$

**Feynman–Kac formula for the backward linearized HJB equation.** Suppose $f : \mathbb{R}^d \times [0,T] \to \mathbb{R}$, $q : \mathbb{R}^d \times [0,T] \to \mathbb{R}^d$ are given functions and let us write a general backward linear advection–diffusion–reaction equation, as it appears in step **(ii)**, in the form:

$$-\partial_t u - \nu \Delta u + q \cdot \nabla u - f = 0. \tag{2.6}$$

The Feynman–Kac representation associated with (2.6) in the time interval $[t_k, t_{k+1}]$ is given by

$$u(x,t_k) = \mathbb{E}\left[ u(X^+_{x,t_k}(t_{k+1}), t_{k+1}) + \int_{t_k}^{t_{k+1}} f(X^+_{x,t_k}(\tau), \tau)\, d\tau \middle| X^+_{x,t_k}(t_k) = x \right]. \tag{2.7}$$

Here, $X^+_{x,t_k}(\cdot)$ denotes a stochastic process initiated at position $x$ at time $t_k$, governed by the drift term $-q$ and by a Wiener process $\sqrt{2\nu}dW$ evolving forward in time over the interval $[t_k, t_{k+1}]$.

As before, we consider a semi-discrete SL-type scheme. Applying again a trapezoidal rule in $[t_k, t_{k+1}]$,

$$\int_{t_k}^{t_{k+1}} f(X^+_{x,t_k}(\tau), \tau)\, d\tau \simeq \frac{\Delta t}{2} \left( f(x, t_k) + f(X^+_{x,t_k}(t_{k+1}), t_{k+1}) \right),$$

in (2.7), we get

$$u(x, t_k) \simeq \mathbb{E}\left[ u\big(X^+_{x,t_k}(t_{k+1}), t_{k+1}\big) + \frac{\Delta t}{2} \left( f(x, t_k) + f\big(X^+_{x,t_k}(t_{k+1}), t_{k+1}\big) \right) \middle| X^+_{x,t_k}(t_k) = x \right]. \tag{2.8}$$

**SL schemes for FP and linearized HJB equations.** In order to obtain a semi-discrete semi-Lagrangian approximation, we approximate the expectations in (2.5) and (2.8) using a weak time discretization of the underlying stochastic processes. Since only expectations are involved, it is sufficient to employ weak schemes, in which the Brownian increment is represented by a finite number of increments that can be interpreted as a discrete random walk approximation of the continuous stochastic process. We consider a weak order 2 Runge-Kutta method, leading to the following one-step time discretization of the stochastic processes $X^-_{x,t_{k+1}}$ and $X^+_{x,t_k}$, given by

$$\Psi^-_\ell(x, t_{k+1}) := x - \frac{\Delta t}{2} \left( q(x, t_{k+1}) + q(x - \Delta t q(x, t_{k+1}) + \xi_\ell, t_k) \right) + \xi_\ell, \tag{2.9}$$

$$\Psi^+_\ell(x, t_k) := x + \frac{\Delta t}{2} \left( q(x, t_k) + q(x + \Delta t q(x, t_k) + \xi_\ell, t_{k+1}) \right) + \xi_\ell, \tag{2.10}$$

respectively, for $\ell \in \mathcal{L}$, where $\mathcal{L}$ is a given finite index set. The vectors $\xi_\ell \in \mathbb{R}^d$ and the associated weights $(\omega_\ell)_{\ell\in\mathcal{L}}$ define a quadrature rule for the Gaussian distribution of the Brownian increment and satisfy $\sum_{\ell\in\mathcal{L}} \omega_\ell = 1$ (see, e.g., [42, 34]).

Let $\mathcal{I}_{\Delta t} = \{0, \ldots, N_T\}$, $\mathcal{I}^*_{\Delta t} = \mathcal{I}_{\Delta t} \setminus \{N_T\}$, and set $\mathcal{G}_{\Delta t} = \{t_k \,|\, k \in \mathcal{I}_{\Delta t}\}$. We denote by $\{m_k : \mathbb{R}^d \to \mathbb{R},\, k \in \mathcal{I}^*_{\Delta t}\}$ a time-discrete approximation of the solution $m$ to (2.3), defined on the temporal grid $\mathcal{G}_{\Delta t}$ according to the following scheme, derived from the representation formula (2.5). Given $m_0(x) = m(x, 0)$, define

$$m_{k+1}(x) = \sum_{\ell\in\mathcal{L}} \omega_\ell\, m_k(\Psi^-_\ell(x, t_{k+1})) \exp\left( \frac{\Delta t}{2} \left( r(\Psi^-_\ell(x, t_{k+1}), t_k) + r(x, t_{k+1}) \right) \right). \tag{2.11}$$

We denote by $\{u_k : \mathbb{R}^d \to \mathbb{R},\, k \in \mathcal{I}^*_{\Delta t}\}$ a discrete approximation of the solution $u$ to (2.6), defined on the temporal grid $\mathcal{G}_{\Delta t}$ according to the following scheme, derived from the representation formula (2.8). Given $u_{N_T}(x) = u(x, T)$, define

$$u_k(x) = \sum_{\ell\in\mathcal{L}} \omega_\ell \left[ u_{k+1}(\Psi^+_\ell(x, t_k)) + \frac{\Delta t}{2} \left( f(\Psi^+_\ell(x, t_k), t_{k+1}) + f(x, t_k) \right) \right]. \tag{2.12}$$

The index set $\mathcal{L}$, the quadrature nodes $(\xi_\ell)_{\ell\in\mathcal{L}}$, and the weights $(\omega_\ell)_{\ell\in\mathcal{L}}$ will be specified in the next sections, in order to construct first and second order schemes.

*Remark* 2.1. In the deterministic case, corresponding to $\nu = 0$, the stochastic perturbations vanish and the expectation in (2.5) and (2.8) reduces to a single contribution. Hence, the quadrature sum disappears and the semi-Lagrangian scheme simply follows the deterministic flow of the drift field $q$. In this case, the Crank–Nicolson one-step maps reduce to

$$\Psi^-(x, t_{k+1}) = x - \tfrac{\Delta t}{2}\Big(q(x, t_{k+1}) + q(x - \Delta t\, q(x, t_k), t_k)\Big),$$
$$\Psi^+(x, t_k) = x + \tfrac{\Delta t}{2}\Big(q(x, t_k) + q(x + \Delta t\, q(x, t_k), t_{k+1})\Big),$$

and the corresponding SL updates become

$$m_{k+1}(x) = m_k(\Psi^-(x, t_{k+1}))\, \exp\Big(\tfrac{\Delta t}{2}\big(r(\Psi^-(x, t_{k+1}), t_k) + r(x, t_{k+1})\big)\Big),$$

$$u_k(x) = u_{k+1}(\Psi^+(x, t_k)) + \tfrac{\Delta t}{2}\Big(f(\Psi^+(x, t_k), t_{k+1}) + f(x, t_k)\Big).$$

These formulas correspond to the standard semi-Lagrangian updates along deterministic characteristics, without any averaging over stochastic realizations. An analogous first order method can be obtained by replacing the Crank–Nicolson discretization with an explicit Euler step in the definition of the flow maps and by employing first order quadrature in time, which leads to simpler but less accurate updates.

## 3. Approximation of the expectation.

In the preceding section, we replaced the exact expectation operator with a quadrature rule. In this section, we present a systematic approach for selecting the quadrature nodes $\{\xi_\ell\}_{\ell\in\mathcal{L}}$ and corresponding weights $\{\omega_\ell\}_{\ell\in\mathcal{L}}$ to achieve both first order and second order accuracy in time.

To ensure clarity and facilitate the derivation, we focus on the prototypical example of the advection-diffusion equation:

$$\partial_t m - \nu\Delta m + q(x)\cdot\nabla m = 0, \tag{3.1}$$

where $\nu > 0$ and, for simplicity, we suppose the drift depending only on the state variable, *i.e.* $q : \mathbb{R}^d \to \mathbb{R}^d$. The same methodology can be extended, with appropriate modifications, to more general advection–diffusion–reaction systems.

Under a semi-Lagrangian discretization, the update (2.11) reduces to

$$m_{k+1}(x) = \sum_{\ell\in\mathcal{L}} \omega_\ell\, m_k(\Psi^-_\ell(x)), \tag{3.2}$$

where the backward flow map is given by

$$\Psi^-_\ell(x) := x - \frac{\Delta t}{2}\left(q(x) + q(x - \Delta t q(x) + \xi_\ell)\right) + \xi_\ell. \tag{3.3}$$

For this test case, we discuss quadrature rules known in the literature and we will derive a new one, to the best of our knowledge, in the framework of second order approximation.

**First order.** To achieve first order accuracy in time, we choose $2d$ axial nodes as

$$\xi_i = \sqrt{2d\nu\Delta t}\, e_i, \quad \xi_{d+i} = -\sqrt{2d\nu\Delta t}\, e_i, \quad i = 1, \dots, d, \tag{3.4}$$

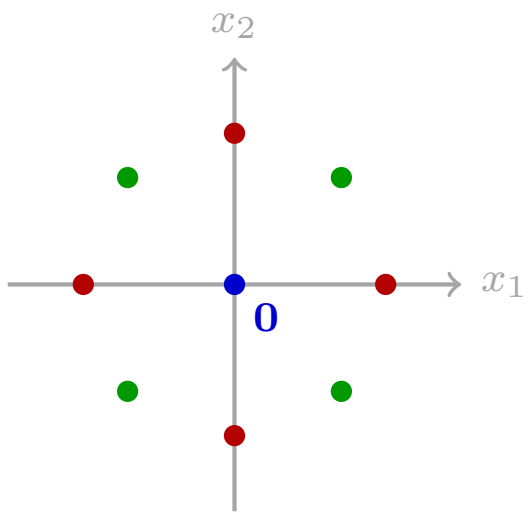

Fig. 1: Illustration for the quadrature nodes for the polynomial rule in dimension 2.

where $\{e_i\}_{i=1}^d$ denotes the canonical basis of $\mathbb{R}^d$. The corresponding weights are chosen as $\omega_\ell = \frac{1}{2d}$ for $\ell = 1, \ldots, 2d$. This method will be denoted as the first order SL scheme (SL1). This quadrature rule exactly matches the moments up to second order of a $d$-dimensional Gaussian distribution with zero mean and covariance $2\nu\Delta t I_{d\times d}$ with probability density $G(y; \Delta t)$. This ensures first order accuracy in time for the scheme (3.2). See, for instance, [8, 17, 41] for a detailed discussion.

**Second order.** To attain second order accuracy in time, it is necessary that the quadrature rule used to approximate Gaussian integrals reproduces the moments of the Gaussian distribution up to the fifth order. This requirement is classical in the analysis of weak approximation schemes for stochastic differential equations; see, for example, Theorem 8.1 in [41].

A natural construction begins with the three-point Gauss-Hermite rule in one dimension, whose nodes reproduce moments up to order five. The $d$-dimensional quadrature nodes $\xi_\ell$ are formed via the full tensor product of the one-dimensional nodes, indexed by $\ell \in \mathcal{L} = \{1, \ldots, 3^d\}$. The corresponding weights $\omega_\ell$ are given by the product of the weights corresponding to each component of $\xi_\ell$.

This scheme will be referred to as the second order exponential rule (SL2e). By construction, this rule exactly reproduces all mixed moments up to the fifth order, yielding a local truncation error of $O(\Delta t^3)$. Its drawback is that the total number of nodes grows as $3^d$, which quickly becomes infeasible for large $d$.

To avoid exponential complexity, one can instead employ a set of nodes whose cardinality grows only quadratically with the dimension $d$. In particular, a symmetric configuration can be adopted, consisting of a central node, axial nodes aligned with the coordinate axes, and face-diagonal nodes corresponding to all pairwise combinations of axes. Specifically, the quadrature rule includes:

(i) A central node $\xi_0 = 0$, with associated weight $\omega_0$;
(ii) $2d$ axial nodes $\xi_i = s\, e_i$ and $\xi_{d+i} = -s\, e_i$ for $i = 1, \ldots, d$, each assigned the same weight $\omega_A$;
(iii) $2d(d-1)$ face-diagonal nodes of the form $\pm s(e_i \pm e_j)$ for all pairs $1 \le i < j \le d$, each assigned weight $\omega_D$,

where $s > 0$ denotes the length of the displacement, which will be specified later along with the weights $\omega_0$, $\omega_A$, and $\omega_D$. A graphical illustration in dimension 2 is given in Figure 1.

This configuration yields a set of nodes indexed by $\mathcal{L} = \{1, \ldots, 2d^2 + 1\}$, since $1 + 2d + 2d(d-1) = 2d^2 + 1$. We refer to this scheme as the second order polynomial rule (SL2p).

The quadrature weights and displacement parameter $s$ are determined so that the rule exactly reproduces the moments of a $d$-dimensional isotropic Gaussian $G(y;\Delta t)$ with zero mean and covariance $2\nu\Delta t\, I_d$, whose moments up to order four are as follows:

$$\text{Zeroth moment (normalization):} \quad \int_{\mathbb{R}^d} G(y;\Delta t)\,dy = 1,$$

$$\text{Second moments (pure):} \quad \int_{\mathbb{R}^d} y_i^2\, G(y;\Delta t)\,dy = 2\nu\Delta t, \qquad i = 1,\dots,d,$$

$$\text{Fourth moments (pure):} \quad \int_{\mathbb{R}^d} y_i^4\, G(y;\Delta t)\,dy = 12(\nu\Delta t)^2, \quad i = 1,\dots,d,$$

$$\text{Fourth moments (mixed):} \quad \int_{\mathbb{R}^d} y_i^2 y_j^2\, G(y;\Delta t)\,dy = 4(\nu\Delta t)^2, \quad i \neq j.$$

Given a multi-index $\alpha$, the following moment conditions

$$\sum_{\ell\in\mathcal{L}} \omega_\ell\, \xi_\ell^\alpha = \int_{\mathbb{R}^d} y^\alpha\, G(y;\Delta t)\,dy, \quad \text{for all } |\alpha| = 0,\dots,5, \tag{3.5}$$

are sufficient to achieve second order consistency, where $|\alpha|$ denotes the order of the multi-index. This is formalized in the following result, with a proof provided in Appendix A.

*Proposition* 3.1. Let $T > 0$, let $m$ be a smooth solution of (3.1) in $\mathbb{R}^d \times [0,T]$ with $q \in C^2(\mathbb{R}^d;\mathbb{R}^d)$, let $\Delta t > 0$ and $k \in \mathcal{I}^*_{\Delta t}$. Given the scheme (3.2), define the corresponding one-step error $\mathcal{T}_k(x) := m(x,t_{k+1}) - \sum_{\ell\in\mathcal{L}} \omega_\ell\, m(\Psi_\ell^-(x,t_{k+1}),t_k)$. Assume the nodes $\{\xi_\ell\}_{\ell\in\mathcal{L}}$ and the weights $\{w_\ell\}_{\ell\in\mathcal{L}}$ satisfy the moment-matching conditions (3.5). Then, there exists a constant $C$, depending only on $T$ and the uniform norms of the derivatives of $m$ and $q$, such that

$$\sup_{x\in\mathbb{R}^d,\, k\in\mathcal{I}^*_{\Delta t}} |\mathcal{T}_k(x)| \le C(\Delta t)^3.$$

The moment-matching conditions (3.5) give rise to a nonlinear algebraic system for the quadrature parameters. Due to the symmetry of the node set, all odd order moments vanish, thereby automatically satisfying (3.5), since the corresponding moments of the Gaussian $G(y;\Delta t)$ are zero. The next proposition records the solution of that system; the detailed computation is deferred to Appendix B.

*Proposition* 3.2. Let $s,\nu,\Delta t > 0$ and $G(y;\Delta t)$ be the $d$-dimensional isotropic Gaussian with zero mean and covariance $2\nu\Delta t\, I_d$. Consider the symmetric quadrature rule with nodes

$$\{0\} \cup \{\pm s e_i\}_{i=1}^d \cup \{\pm s(e_i \pm e_j)\}_{1\le i<j\le d},$$

and associated weights $\omega_0$ (for the central node), $\omega_A$ (for each axial node $\pm s e_i$) and $\omega_D$ (for each face-diagonal node $\pm s(e_i \pm e_j)$). If the parameters $s$, $\omega_0$, $\omega_A$, $\omega_D$ are chosen as

$$s = \sqrt{6\nu\Delta t}, \quad \omega_0 = \frac{d^2 - 7d + 18}{18}, \quad \omega_A = \frac{4-d}{18}, \quad \omega_D = \frac{1}{36},$$

then the corresponding quadrature rule satisfies the moment-matching conditions (3.5).

*Remark* 3.3. Note that $\omega_A = (4-d)/18$ is negative for $d > 4$, so for these dimensions the axial weights do not admit a probabilistic interpretation; this may have consequences for positivity and stability. Nevertheless, numerical experiments indicate that the scheme is asymptotically positivity-preserving, meaning that violations of positivity vanish as $\Delta t \to 0$. The implications of negative weights for stability can be addressed by standard remedies (e.g., alternative node sets or constrained weight adjustments) if required.

## 4. Tensor-Train Approximation.

In the previous sections, we derived a semi-Lagrangian scheme for approximating the MFG system. This method requires evaluating the approximate solution at previous or next time steps: at $\Psi_\ell^-(\cdot,\cdot)$ for the FP equation, and at $\Psi_\ell^+(\cdot,\cdot)$ for the HJB equation. In the classical SL framework, these evaluations are obtained via interpolation, typically piecewise linear, based on grid point values. Although efficient in low dimensions, this approach becomes infeasible in high dimensions, which is the main focus of this work.

To overcome this challenge, we consider a Tensor-Train approximation to efficiently represent high-dimensional functions. The remainder of this section introduces the TT methodology in three steps:

- We first describe how functions are represented in a multivariate Legendre polynomial basis.
- Next, we present the TT representation, which expresses the coefficient tensor in a low-rank factorized form.
- Finally, we summarize practical TT construction algorithms that enable efficient approximation from sampled data.

Together, these subsections provide the necessary background for the Tensor-Train-based semi-Lagrangian method developed in Section 5.

**Polynomial Basis Representation.** Let $\Pi_{i,n_i} = \operatorname{span}\{\psi_{i,1},\dots,\psi_{i,n_i}\}$, for $i = 1,\dots,d$, be the space of univariate Legendre polynomials of degree at most $n_i$, where $\psi_{i,k}(\cdot)$ denotes the $k$-th Legendre polynomial in the $i$-th variable. Their tensor product defines the multivariate polynomial space

$$\mathcal{V}_d = \Pi_{1,n_1} \otimes \cdots \otimes \Pi_{d,n_d}.$$

A general function $v \in \mathcal{V}_d$ can then be expressed as

$$v(x_1,\dots,x_d) = \sum_{i_1=1}^{n_1} \cdots \sum_{i_d=1}^{n_d} c[i_1,\dots,i_d]\, \psi_{1,i_1}(x_1)\cdots\psi_{d,i_d}(x_d), \tag{4.1}$$

where $c \in \mathbb{R}^{n_1\times\cdots\times n_d}$ is the multidimensional coefficient tensor.

Legendre polynomials form an orthogonal basis of $L^2([-1,1])$ and provide efficient approximations for smooth functions on bounded intervals. Let $n = \max\{n_k\}$ denote the maximum polynomial degree used in the representation. For functions with $s$ continuous derivatives, the Legendre projection error decays as $\mathcal{O}(n^{-s})$, while for analytic functions the convergence is exponential in $n$ [12]. These properties make Legendre expansions particularly effective for representing smooth functions in high-dimensional settings.

However, the size of $c$ grows exponentially with the dimension $d$, leading to the curse of dimensionality.

**Tensor-Train Representation.** The TT format mitigates this problem by factorizing $c$ into a sequence of lower order tensors (called TT cores). For given TT ranks $\{r_k\}_{k=1}^{d-1}$, we define:

$$U_1 \in \mathbb{R}^{n_1 \times r_1}, \quad U_k \in \mathbb{R}^{r_{k-1} \times n_k \times r_k} \ (k = 2, \dots, d-1), \quad U_d \in \mathbb{R}^{r_{d-1} \times n_d}.$$

The entries of $c$ are approximated as

$$\tilde{c}[i_1, \dots, i_d] = \sum_{j_1=1}^{r_1} \cdots \sum_{j_{d-1}=1}^{r_{d-1}} U_1[i_1, j_1]\, U_2[j_1, i_2, j_2] \cdots U_d[j_{d-1}, i_d]. \tag{4.2}$$

Accordingly, the TT approximation of $v$ reads:

$$\tilde{v}(x_1, \dots, x_d) = \sum_{i_1=0}^{n_1} \cdots \sum_{i_d=0}^{n_d} \tilde{c}[i_1, \dots, i_d]\, \psi_{i_1}(x_1) \cdots \psi_{i_d}(x_d). \tag{4.3}$$

If we let $R = \max\{r_k\}$, the number of parameters in this representation scales as $\mathcal{O}(dnR^2)$, linearly in $d$ for fixed ranks. This efficient scaling makes the TT format suitable for high-dimensional problems.

Many important function classes, such as quadratic functions [23, 40] and weakly correlated Gaussian functions [51], admit low-rank TT representations. This property is particularly useful in MFG applications, where quadratic value functions and Gaussian densities naturally arise. In the Linear–Quadratic–Gaussian MFG setting, characterized by linear stochastic dynamics and quadratic cost structures, the HJB equation admits quadratic solutions in the state variable, while the associated FP equation preserves the Gaussian structure of the population density, making this framework widely used as a benchmark model in the literature [31, 32].

**TT Construction Methods.** The TT approximation can be obtained via algorithms based on the Singular Value Decomposition (SVD). The TT-SVD algorithm [47] constructs a TT decomposition with guaranteed accuracy but requires full access to all tensor entries, an operation that quickly becomes intractable in high dimensions.

To avoid this bottleneck, we employ the TT-Cross approximation method. Here, the TT representation is determined from function evaluations at a finite set of sampling points $\{z_p\}_{p=1}^N \subset \mathbb{R}^d$ by solving the least-squares problem:

$$\min_{U_1, \dots, U_d} \sum_{p=1}^{N} |\tilde{v}(z_p) - v(z_p)|^2 , \tag{4.4}$$

where $\tilde{v}$ depends on the TT cores $\{U_i\}_{i=1}^d$ through (4.2)–(4.3).

The TT cores are updated alternately over the dimensions using cross interpolation techniques [48], which require only function samples rather than the full tensor grid. Pivoting strategies, such as the `maxvol` algorithm [30], are used to select well-conditioned sample sets, referred to as adaptive sampling sets, that improve stability and accuracy. This alternating iteration continues until convergence, typically measured by the relative change in the TT cores. A detailed overview is provided in [24], and the quasi-optimality of the method is established in [53].

Finally, once in TT form, function values, derivatives, and integrals can be efficiently computed. In particular, gradients $\nabla \tilde{v}(x)$ can be evaluated at a cost $\mathcal{O}(dnR^2)$ per point, an essential feature for high-dimensional control and optimization. Further details on TT gradient evaluation can be found in [49].

**5. Tensor-Train-Based Semi-Lagrangian Method.** We now describe a numerical scheme for high-dimensional MFGs that combines a semi-Lagrangian time discretization with TT function representations. This approach enables efficient implementation of the policy iteration algorithm SPI in high dimensions. Both the density $m$ and the value function $u$ are represented in a Legendre polynomial tensor basis, with their coefficient tensors stored in TT format. This provides compact storage and allows efficient manipulation of high-dimensional data. The semi-Lagrangian discretization updates the functions in time, while the TT format ensures scalable low-rank representations.

We first introduce the notation and then summarize the computational steps in Algorithm SPISL.

**Notation.** For $(n, x, k, \ell) \in \mathbb{N} \times \mathbb{R}^d \times \mathcal{I}^*_{\Delta t} \times \mathcal{L}$, we denote by $\Psi^{(n),+}_{\ell,k}(x)$ and $\Psi^{(n),-}_{\ell,k+1}(x)$ the forward and backward maps at iteration $n$, corresponding to drifts $\overline{q}^{(n)}_k$ and $\overline{q}^{(n)}_{k+1}$. The running cost is defined as

$$f^{(n)}_k(x) = L\big(\overline{q}^{(n)}_k(x)\big) + F\big(x, m^{(n)}_k\big).$$

**Algorithm SPISL.** Let $q^{(0)} : \mathbb{R}^d \times \mathcal{G}_{\Delta t} \to \mathbb{R}^d$ be an initial policy and set $\overline{q}^{(0)} = q^{(0)}$. For each iteration $n \geq 0$, given $\delta_n \in (0, 1)$, one outer policy iteration cycle consists of:

1. **Forward (FP) update.** For each $k \in \mathcal{I}^*_{\Delta t}$, define

$$m^{(n)}_{k+1}(x) = \sum_{\ell \in \mathcal{L}} \omega_\ell\, \widetilde{m}^{(n)}_k\big(\Psi^{(n),-}_{\ell,k+1}(x)\big) \exp\big(\Delta t\, r^{(n)}_{\ell,k}(x)\big), \qquad x \in \mathbb{R}^d, \tag{5.1}$$

   with

$$r^{(n)}_{\ell,k}(x) = \tfrac{1}{2}\Big(\mathrm{div}(\overline{q}^{(n)}_k(\Psi^{(n),-}_{\ell,k+1}(x))) + \mathrm{div}(\overline{q}^{(n)}_{k+1}(x))\Big), \quad \widetilde{m}^{(n)}_0 = m_0.$$

   From the sample pairs $(x, m^{(n)}_{k+1}(x))$, we solve the least-squares problem (4.4) to obtain TT cores and reconstruct $\widetilde{m}^{(n)}_{k+1}$. Off-grid evaluations of $\widetilde{m}^{(n)}_k(\cdot)$ use TT point evaluation (4.3).
2. **Backward (HJB) update.** For each $k \in \mathcal{I}^*_{\Delta t}$, define

$$u^{(n)}_k(x) = \textstyle\sum_{\ell \in \mathcal{L}} \omega_\ell \Big(\widetilde{u}^{(n)}_{k+1}(\Psi^{(n),+}_{\ell,k}(x)) + \frac{\Delta t}{2} f^{(n)}_{k+1}(\Psi^{(n),+}_{\ell,k}(x))\Big) + \frac{\Delta t}{2} f^{(n)}_k(x), \qquad x \in \mathbb{R}^d, \tag{5.2}$$

   with terminal condition $\widetilde{u}^{(n)}_{N_T}(x) = G(x, \widetilde{m}^{(n)}_{N_T})$. Least-squares reconstruction and TT evaluation are performed analogously to the FP step.
3. **Policy update.** For each $k \in \mathcal{I}_{\Delta t}$, the control policy is updated as

$$q^{(n+1)}_k(x) = \nabla_p H\big(\nabla \widetilde{u}^{(n)}_k(x)\big), \quad x \in \mathbb{R}^d.$$

   Gradients $\nabla \widetilde{u}^{(n)}_k(x)$ are efficiently computed in TT form at a cost $\mathcal{O}(dnR^2)$ per point.
4. **Smoothing.** To stabilize iterations, apply relaxation:

$$\overline{q}^{(n+1)}_k(x) = (1 - \delta_n)\, \overline{q}^{(n)}_k(x) + \delta_n\, q^{(n+1)}_k(x).$$

**Implementation details.** In practice, several aspects of the algorithm can be handled efficiently to improve stability and convergence:

- *Initial policy.* A common initialization for the control policy is
$$q_k^{(0)}(x) = \nabla_p H\big(\nabla G(x)\big), \quad k \in \mathcal{I}_{\Delta t},$$
that is, by using the gradient of the terminal Hamiltonian condition. If this terminal condition depends on the FP solution, it can be evaluated at the initial FP state, exploiting the known boundary data.
- *Stopping criterion.* Convergence of the outer policy iteration is monitored using a fixed tolerance $\varepsilon_{\text{stop}} > 0$. At the end of each iteration, the algorithm checks whether the sum of the $L^2$-norm differences between two consecutive iterates, the FP density at the final time and the HJB solution at the initial time, falls below this tolerance, namely
$$\big\|m_{N_T}^{(n)} - m_{N_T}^{(n-1)}\big\|_{L^2(\mathbb{R}^d)} + \big\|u_0^{(n)} - u_0^{(n-1)}\big\|_{L^2(\mathbb{R}^d)} \leq \varepsilon_{\text{stop}}.$$
This condition ensures that both solution components have reached a consistent fixed point at their respective time boundaries.
- *TT-Cross initialization.* The TT-Cross reconstructions at each iteration can be initialized using the TT cores and sample sets obtained in the previous policy iteration. This warm-start strategy significantly accelerates the convergence of the least-squares reconstruction steps.
- *Logarithmic transform.* For Gaussian-like densities, it is often advantageous to perform a logarithmic transform before TT reconstruction. Since $\log m$ of a Gaussian is a quadratic function, it can be represented compactly in the Legendre basis. Given samples of $m$, one computes $m_{\log} = \log m$, reconstructs $m_{\log}$ in TT form by solving the least-squares problem (4.4), and recovers $m(x) = \exp(m_{\log}(x))$ when needed. This strategy improves compressibility and ensures positivity of the density by construction.

**6. Numerical tests.** We evaluate the efficiency and accuracy of the proposed Tensor-Train-based semi-Lagrangian scheme across a representative suite of benchmark problems. We consider both local and nonlocal mean field game problems, systematically varying spatial dimension, temporal resolution, and quadrature strategy. Our analysis focuses on numerical error in both $L^2$ and uniform norms, empirical convergence rates, and computational cost measured in CPU time. The space-time domain is denoted by $Q = [-L, L]^d \times [0, T]$, for a given constant $L > 0$. Numerical boundary conditions are handled via extrapolation; see [16] for further details. The numerical errors are evaluated as follows. Let $\{x_k\}_{k=1}^{N_c}$ be a set of $N_c = 10^5$ validation points uniformly distributed in the domain $[-L, L]^d$. We define the discrete $L^2$ and uniform errors as

$$E_2(h) = \left( \frac{\sum_{k=1}^{N_c} \big(h(x_k) - h_{\text{ex}}(x_k)\big)^2}{\sum_{k=1}^{N_c} \big(h_{\text{ex}}(x_k)\big)^2} \right)^{1/2}, \; E_\infty(h) = \max_{1 \leq k \leq N_c} \frac{\big|h(x_k) - h_{\text{ex}}(x_k)\big|}{\big|h_{\text{ex}}(x_k)\big|},$$

where $h$ denotes the numerical solution, and $h_{\text{ex}}$ denotes the exact solution. When a forward equation for the density is considered, the errors $E_2(m)$ and $E_\infty(m)$ are evaluated at the final time $t = T$. In the coupled Mean Field Game experiments, the errors for the value function $u$ are computed at the initial time $t = 0$, consistently with the backward–forward structure of the system. These indicators will be used to compute numerical convergence rates and to verify whether it matches the expected values. Note that the numerical order will be evaluated only with respect to time, since the number of spatial reconstruction polynomials is fixed for all the numerical tests.

In Table 1, we list the schemes analyzed in the following sections for the viscous-case problems, along with their number of characteristics. Here, SL1 denotes the first order SL scheme with $2d$ characteristics, SL2e denotes the standard second order SL scheme with $3^d$ characteristics, and SL2p denotes the new second order SL scheme with a polynomial number of characteristics.

| **Scheme** | **# Characteristics** |
|---|---|
| SL1 - First Order | $2d$ |
| SL2e - Second Order (Exponential) | $3^d$ |
| SL2p - Second Order (Polynomial) | $2d^2+1$ |

Table 1: Semi-Lagrangian schemes considered and the corresponding number of characteristics.

We denote by $n_{\text{HJB}}$ and $n_{\text{FP}}$ the dimensions of the Legendre expansions employed in the TT approximations of the solutions to the HJB and FP equations, respectively. These dimensions are selected such that the relative $L^2$–error of the reconstructed initial condition for the FP equation and the final condition for the HJB equation, evaluated over the $N_c$ validation nodes, remains below a prescribed tolerance $tol_{start}$.

All the tests have been performed on a Dell XPS 13 with an Intel Core i7 (2.8 GHz) and 16 GB RAM. The codes are written in MATLAB R2023b.

**6.1. Advection-Diffusion Equation with Exact Solution.** We begin by considering a simple test case, adapted from [28], based on the linear advection-diffusion equation:

$$\partial_t m - \nu\Delta m + q\cdot\nabla m = 0,$$

subject to periodic boundary conditions, with $q\in\mathbb{R}^d$ a constant advection field. The exact solution is given by

$$m(x,t) = 2 + \sin\left(\pi\sum_{i=1}^{d}(x_i - t - s_i)\right)e^{-\nu\pi^2 d\, t},$$

where $s_i = \frac{i-1}{d}$ for $i = 1,\dots,d$.

By trigonometric identities, the rank of this solution can be shown to be $2^{d-1}$, which rapidly increases with the spatial dimension $d$. We fix $tol_{start} = 10^{-8}$, yielding $n_{\text{FP}} = 15$ basis functions for the TT reconstruction, $Q = [-1,1]^d\times[0,T]$, $q = (1,\dots,1)$ and $\nu = 0.1$. As in [28], we define the final time as $T = \log(2)/(d\nu\pi^2)$, which represents the time at which the diffusion process reduces the amplitude of the solution by half. In Table 2, we show, for several refinements of the time mesh, the CPU time, the error $E_2(m)$, and the corresponding numerical convergence rate for the second order schemes SL2e and SL2p in dimension $d = 3$.

We observe that the SL2e consistently yields lower error and marginally higher computational cost for coarse $\Delta t$, whereas the SL2p becomes competitive as $\Delta t$ decreases. In both cases, the numerical convergence rates remain very close to the expected value of two. We note that the time step values reported in Table 2 correspond to a refinement of the time mesh starting from just 2 steps up to 32. Even the coarsest choice (2 time steps) provides accurate results due to the semi-Lagrangian nature of the schemes, which are not restricted by CFL conditions.

| | SL2e | | | SL2p | | |
|---|---|---|---|---|---|---|
| $\Delta t$ | CPU (s) | $E_2(m)$ | Order | CPU (s) | $E_2(m)$ | Order |
| 0.1171 | 1.25 | $1.25 \times 10^{-4}$ | – | 0.91 | $8.23 \times 10^{-4}$ | – |
| 0.0585 | 2.00 | $3.02 \times 10^{-5}$ | 2.05 | 1.58 | $1.90 \times 10^{-4}$ | 2.11 |
| 0.0293 | 3.53 | $7.43 \times 10^{-6}$ | 2.02 | 2.75 | $4.57 \times 10^{-5}$ | 2.06 |
| 0.0146 | 7.04 | $1.84 \times 10^{-6}$ | 2.01 | 5.11 | $1.12 \times 10^{-5}$ | 2.03 |
| 0.0073 | 15.51 | $4.59 \times 10^{-7}$ | 2.01 | 10.16 | $2.78 \times 10^{-6}$ | 2.01 |

Table 2: Test 6.1. Comparison of the SL2e vs. SL2p schemes for the advection–diffusion equation in the case $\nu = 0.1$, $d = 3$, with $n_{\mathrm{FP}} = 15$.

We now assess the schemes' performance as the spatial dimension $d$ increases, while keeping the time-step size $\Delta t$ fixed. Specifically, for the first order scheme SL1 we set $\Delta t = T/128$, whereas for the second order methods we choose $\Delta t = T/2$ for the exponential SL2e scheme and $\Delta t = T/4$ for the polynomial SL2p scheme. This choice ensures that all schemes are compared under time discretizations that yield errors of comparable magnitude, namely around $3 \cdot 10^{-4}$.

| $d$ | SL1 | SL2e | SL2p |
|---|---|---|---|
| 3 | 14.21 | **0.70** | 1.02 |
| 4 | 35.65 | 3.80 | **3.03** |
| 5 | 99.89 | 19.10 | **7.17** |
| 6 | 175.40 | 80.30 | **16.43** |
| 7 | 219.76 | 370.00 | **30.05** |
| 8 | 223.30 | 1412.50 | **51.78** |

Table 3: Test 6.1. CPU times (in seconds) for the SL1, SL2e, and SL2p schemes in the case $\nu = 0.1$, varying the dimension $d$.

In Table 3 we show the CPU times (in seconds) and the lowest value in every row is highlighted in bold. The SL2p scheme becomes the most efficient method already from $d = 4$ onward, consistently outperforming both SL1 and SL2e as the dimension increases. Consequently, for low dimensions the scheme SL2e provides a favorable trade-off between accuracy and runtime, whereas for higher dimensions ($d > 3$) the scheme SL2p becomes preferable. In all cases, the first-order scheme SL1 is strictly outperformed by its second-order counterparts in both accuracy and efficiency. Remarkably, in all test cases and for all time-step sizes, the TT approximation maintained a uniformly low rank, with a maximum TT rank equal to 4, confirming the high compressibility of the solutions considered.

We analyze the computational performance of each method by fitting the CPU time $C(d)$ as a function of the problem dimension $d$ using two models. The first is an exponential model, $C(d) = a\, e^{b\, d}$, which is estimated via a semi-logarithmic regression of $C(d)$ versus $d$. The second is a power-law model, $C(d) = a\, d^b$, obtained through a log-log regression of $C(d)$ versus $d$. To evaluate the accuracy of the fitted models, we compute the coefficient of determination ($R^2$), defined as $R^2 = 1 - \frac{\text{RSS}}{\text{TSS}}$. Here, RSS and TSS denote the Residual Sum of Squares and Total Sum of Squares, respectively. For a set of observed values $\{y_i\}$ and predicted values $\{\hat{y}_i\}$, we have

$$\text{RSS} = \sum_i (y_i - \hat{y}_i)^2, \qquad \text{TSS} = \sum_i (y_i - \bar{y})^2,$$

where $\bar{y}$ is the mean of the observed values. Values of $R^2$ close to 1 indicate that the model explains most of the variability in the data, while values near 0 indicate a poor fit.

| Method | Fitting | Coefficient $a$ | Exponent $b$ | $R^2$ |
|---|---|---|---|---|
| SL1 | Exponential | 2.00 | 0.61 | 0.963 |
| SL1 | Polynomial | 0.31 | 3.18 | **0.994** |
| SL2e | Exponential | 0.008 | 1.52 | **0.976** |
| SL2e | Polynomial | 9.91e-05 | 7.75 | 0.881 |
| SL2p | Exponential | 0.12 | 0.78 | 0.915 |
| SL2p | Polynomial | 0.012 | 4.03 | **0.999** |

Table 4: Test 6.1. Results of the fittings for each scheme using both exponential and power-law models for the case $\nu = 0.1$.

Table 4 shows, for each scheme, the fitted parameters $a$ and $b$ for both models, along with the corresponding coefficient of determination. For each scheme, the highest $R^2$ is highlighted in bold. For the first order scheme, the polynomial model yields the most accurate fit with $C(d) \approx 0.31\, d^{3.18}$, corresponding to a near-cubic scaling of computational cost. The SL2e scheme displays the opposite trend: the exponential fit performs significantly better, $C(d) \approx 0.008\, e^{1.52\, d}$, indicating a strong exponential dependence of the computational cost on the dimension. Finally, for the SL2p scheme, the polynomial fit provides an excellent representation of the data, $C(d) \approx 0.012\, d^{4.03}$, suggesting a quartic scaling $O(d^4)$.

Overall, these results confirm that the exponential model best describes the SL2e behavior, while polynomial growth accurately characterizes both SL1 and SL2p schemes. Among the second order methods, SL2p demonstrates a superior balance between accuracy and computational scaling up to moderate dimensions, whereas SL1 may become more efficient in very high-dimensional regimes due to its slower polynomial growth.

**Positivity test.** In this section, we examine the preservation of the positivity of the solution. As discussed in Remark 3.3, the weight $\omega_A$ becomes negative for spatial dimensions $d > 4$, potentially compromising both stability and positivity. Furthermore, the magnitude of this negativity increases with the dimension $d$. Consequently, we focus our analysis on the largest dimension considered in the previous test case, namely $d = 8$.

We consider the following exact solution in $[-1, 1]^8 \times [0, T]$:

$$m(x,t) = 0.5 + \sin\left(\pi \sum_{i=1}^{d} (x_i - t)\right) \exp\bigl(-d\nu\pi^2 t\bigr).$$

By construction, the minimum of this solution at the final time $T = \ln(2)/(d\nu\pi^2)$ is zero, which occurs, for example, at the points $x_k^* = (s_k, \ldots, s_k)$, where $s_k = T +$ $k/4 - 1/16$ and $k \in I = \{-4, \dots, 3\}$. Throughout this test, we employ the same parameter values as in the previous case.

To assess the preservation of positivity, we evaluate the numerical approximation at the final time, denoted by $m_{N_t}$, at the points $\{x_k^*\}_{k\in I}$ and record the minimum value. The results are summarized in Table 5. As anticipated, the numerical solution exhibits a slight loss of positivity; however, the minimum value decreases with a rate close to 2 as the time step $\Delta t$ is refined. This behavior is consistent with the theoretical second order accuracy of the SL2p scheme and indicates that the method is asymptotically positivity-preserving, since positivity is recovered in the limit $\Delta t \to 0$.

| $\Delta t$ | $\min_k m_{N_t}(x_k^*)$ | Order |
|---|---|---|
| 0.0439 | -5.78e-03 | |
| 0.0219 | -1.29e-03 | 2.17 |
| 0.0110 | -3.04e-04 | 2.08 |

Table 5: Test 6.1. Positivity test analysis for the SL2p scheme varying time step $\Delta t$ in the case $\nu = 0.1$ and $d = 8$.

**6.2. Local MFG with exact solution.** We consider the following local MFG system:

$$\begin{cases} -\partial_t u - \nu\Delta u + \frac{|\nabla u|^2}{2} - \frac{\beta |x|^2}{2} & = \gamma \ln m, \\ \partial_t m - \nu \Delta m - \operatorname{div}(m\nabla u) & = 0, \end{cases} \tag{6.1}$$

subject to the terminal and initial conditions:

$$u(x,T) = \frac{\alpha|x|^2}{2} - \left(\nu d\alpha + \frac{\gamma d}{2}\ln\frac{\alpha}{2\pi\nu}\right)T, \quad m(x,0) = \left(\frac{\alpha}{2\pi\nu}\right)^{\frac{d}{2}} e^{-\frac{\alpha|x|^2}{2\nu}},$$

where $\nu > 0$ and $\alpha = \frac{-\gamma + \sqrt{\gamma^2 + 4\nu^2\beta}}{2\nu}$. In [39], the authors demonstrate that such system admits the explicit solution

$$u(x,t) = \frac{\alpha|x|^2}{2} - \left(\nu d\alpha + \frac{\gamma d}{2}\ln\frac{\alpha}{2\pi\nu}\right)t, \quad m(x,t) = \left(\frac{\alpha}{2\pi\nu}\right)^{\frac{d}{2}} e^{-\frac{\alpha|x|^2}{2\nu}},$$

where $m$ remains stationary and $u$ evolves linearly in time.

The parameter $\alpha$ decreases with the congestion parameter $\gamma$, reflecting the effect of the logarithmic term $\gamma \ln m$ in the HJB equation: stronger congestion (larger $\gamma$) leads to a more diffuse density profile.

We begin our analysis with the resolution of the HJB equation, first comparing it with grid-based techniques, and subsequently testing it in a high-dimensional setting. Afterwards, we turn our attention to the coupled MFG system.

**Comparison with Grid-Based SL Scheme.** Let us consider the HJB equation in (6.1) posed on the three-dimensional cubic domain with final time $T = 0.02$, $Q = [-0.1, 0.1]^3 \times [0, 0.02]$, viscosity coefficient $\nu = 0.01$, and parameter $\gamma = 0$.

In this case, where only the HJB equation is approximated, we apply a modified version of Algorithm SPISL. In particular, Step 1 is omitted, while Step 2 is performed by updating, at each time step $k$, the drift appearing in $\Psi_{\ell,k}^{(n),+}(x)$ according to

$$q_k^{(n)}(x) = \nabla_p H\big(\nabla \widetilde{u}_{k+1}^{(n)}(x)\big), \qquad x \in \mathbb{R}^d,$$

where the function $\widetilde{u}_{k+1}^{(n)}$ is the approximation computed at the previous time step. In our numerical experiments, this procedure requires only one outer policy iteration to satisfy the stopping criterion. The objective of this experiment is to evaluate the accuracy and efficiency of the proposed method in comparison with a classical grid-based semi-Lagrangian approach. The choice of this relatively small domain and final time is made to obtain a sufficiently accurate approximation within a reasonable CPU time for the grid-based semi-Lagrangian scheme. Both schemes employ a first order SL discretization. Table 6 presents the $L^2$ and uniform errors of the TT approximation for decreasing values of the time step $\Delta t$, along with the corresponding CPU times. The number of Legendre basis functions used in the functional TT approximation is fixed to three. In this case, the TT ranks remain uniformly equal to 3 for all time steps. For reference, Table 7 reports the errors and CPU times obtained using a classical grid-based semi-Lagrangian method, implemented on a uniform Cartesian grid with $n$ points per spatial dimension. In this case, both the time step and spatial resolution are refined simultaneously.

| $\Delta t$ | $E_2(u)$ | $E_\infty(u)$ | CPU (s) |
|---|---|---|---|
| 0.0056 | $9.33 \times 10^{-5}$ | $9.21 \times 10^{-5}$ | 0.14 |
| 0.0026 | $4.21 \times 10^{-5}$ | $4.14 \times 10^{-5}$ | 0.27 |
| 0.0013 | $2.00 \times 10^{-5}$ | $1.97 \times 10^{-5}$ | 0.52 |

Table 6: Test 6.2, HJB equation. $L^2$, $L^\infty$ errors and CPU time for the SL1 scheme in the case $d = 3$, with $n_{HJB} = 3$.

| $n$ | $\Delta t$ | $E_2(u)$ | $E_\infty(u)$ | CPU (s) |
|---|---|---|---|---|
| 10 | 0.0056 | $9.36 \times 10^{-5}$ | $9.26 \times 10^{-5}$ | 1.08 |
| 20 | 0.0026 | $4.94 \times 10^{-5}$ | $4.85 \times 10^{-5}$ | 37.29 |
| 40 | 0.0013 | $2.52 \times 10^{-5}$ | $2.47 \times 10^{-5}$ | 2057.09 |

Table 7: Test 6.2, HJB equation. $L^2$, $L^\infty$ errors and CPU time for the SL1 grid-based scheme in the case $d = 3$.

Both methods attain first-order convergence in the time-step, as evidenced by the approximately halved errors when $\Delta t$ is halved. Crucially, the TT solver achieves comparable accuracy to the grid-based scheme at a fraction of the computational cost: the CPU time for the TT method grows approximately linearly as $\Delta t$ decreases, whereas the grid solver's CPU time increases superlinearly due to the $O(n^3)$ cost of refining the spatial mesh. In particular, reducing $\Delta t$ from $5.6 \times 10^{-3}$ to $1.3 \times 10^{-3}$ multiplies the TT CPU time by under four, while the grid-based CPU time grows by nearly three orders of magnitude. These results highlight the pronounced efficiency of the TT compression in alleviating the curse of dimensionality, offering substantial computational advantages already for three-dimensional problems and beyond.

**High-dimensional test.** We now investigate the performance of the proposed technique in solving a high-dimensional HJB equation. We fix $T = 1$, $Q = [-1, 1]^d \times [0, 1]$, $\nu = 1$ and $\gamma = 0.1$. For the time discretization, we employ the first order SL1 scheme. This choice is motivated by the fact that, although the second order SL2p method achieves greater accuracy in moderate dimensions, its computational cost scales with a polynomial of higher degree compared to SL1. Consequently, SL2p becomes less practical in very high-dimensional settings. In contrast, SL1 provides a more favorable compromise between accuracy and computational efficiency when $d$ is large.

| $d$ | $E_2(u)$ | $E_\infty(u)$ | CPU (s) |
|---|---|---|---|
| 50 | $2.48 \times 10^{-2}$ | $2.11 \times 10^{-2}$ | $1.63 \times 10^{3}$ |
| 100 | $2.49 \times 10^{-2}$ | $2.16 \times 10^{-2}$ | $1.10 \times 10^{4}$ |

Table 8: Test 6.2, HJB equation. $L^2, L^\infty$ errors and CPU times for the SL1 scheme varying dimension $d$ with $\Delta t = 0.0312$ and $n_{HJB} = 3$.

The numerical results are summarized in Table 8. For $d = 50$, the computation required approximately 30 minutes, whereas the $d = 100$ case took around 3 hours. The accuracy, measured in both $L^2$ and $L^\infty$ norms, remains of order $10^{-2}$ and essentially stable as the dimension doubles from $d = 50$ to $d = 100$. This indicates that the error is not strongly affected by the growth in dimensionality, at least within the tested range, which is a promising feature of the approach. We remark that this level of accuracy is consistent with the first-order approximation used in the scheme; similar errors are observed in lower-dimensional cases (e.g., $d = 3$), confirming that the method order, rather than the dimensionality, dominates the error. In the following section, we demonstrate how a second-order scheme improves the accuracy.

In addition, the TT ranks remain uniformly equal to 3 in all tested cases. Finally, we note that reducing the time step $\Delta t$ leads to smaller errors, as expected, but at the expense of a substantial increase in CPU time. Thus, there is a clear trade-off between accuracy and efficiency.

**Mean Field Games test.** We now turn to the analysis of the full MFG system, coupling the HJB and FP equations. For these tests, we employ the SL2p scheme, which provides second order accuracy in time, and consider the cases $d = 3$ and $d = 6$. We fix $Q = [-1, 1]^d \times [0, T]$, with $T = 1$, $\nu = 1$ and $\gamma = 0.1$. To handle the FP equation, we employ the log–exp–log transformation so that the density remains a quadratic polynomial in the TT representation. Further details on this methodology are provided in Section 5. In this setting, we select the stopping tolerance $\varepsilon_{\text{stop}} = 10^{-5}$ and the smoothing parameter $\delta_n = 10^{-2}$, motivated by numerical evidence indicating that smaller values of $\delta_n$ accelerate the convergence of the policy iteration.

| $\Delta t$ | $E_2(u)$ | Order | $E_2(m)$ | Order | CPU (s) |
|---|---|---|---|---|---|
| $2.50 \times 10^{-1}$ | $1.21 \times 10^{-2}$ | – | $2.00 \times 10^{-3}$ | – | 2.25 |
| $1.25 \times 10^{-1}$ | $2.81 \times 10^{-3}$ | 2.11 | $5.36 \times 10^{-4}$ | 1.90 | 5.25 |
| $6.25 \times 10^{-2}$ | $6.77 \times 10^{-4}$ | 2.05 | $1.35 \times 10^{-4}$ | 1.99 | 3.52 |

Table 9: Test 6.2, local MFG system. $L^2$ errors and CPU time for the SL2p scheme in the case $d = 3$, with $n_{HJB} = n_{FP} = 3$.

Tables 9 and 10 report the $L^2$ errors for both the HJB and FP components, along with the observed convergence orders and CPU times. In both cases, we clearly observe second order convergence, confirming that the SL2p scheme delivers the expected accuracy for the coupled system. The TT rank associated with the approximation of the value function remains fixed at 3 for all time steps and all the iterations, while the maximal TT rank required for the approximation of the FP equation is equal to 9.

| $\Delta t$ | $E_2(u)$ | Order | $E_2(m)$ | Order | CPU (s) |
|---|---|---|---|---|---|
| $2.50 \times 10^{-1}$ | $1.30 \times 10^{-2}$ | – | $4.40 \times 10^{-3}$ | – | 17.33 |
| $1.25 \times 10^{-1}$ | $3.05 \times 10^{-3}$ | 2.10 | $1.10 \times 10^{-3}$ | 1.99 | 17.55 |
| $6.25 \times 10^{-2}$ | $7.31 \times 10^{-4}$ | 2.06 | $2.62 \times 10^{-4}$ | 2.08 | 36.38 |

Table 10: Test 6.2, local MFG system. $L^2$ errors and CPU time for the SL2p scheme in the case $d = 6$, with $n_{HJB} = n_{FP} = 3$.

Regarding the CPU times, we observe at most a linear growth with respect to the number of time steps. Interestingly, in certain cases the cost remains nearly invariant, or even decreases, when the time step is refined. This effect can likely be attributed to the TT-Cross reconstructions, which at each iteration are initialized using the TT elements from the previous step. Such a warm-start strategy accelerates the reconstruction procedure and can explain the observed behavior.

**6.3. Non-local MFG with exact solution.** We consider the following non-local MFG system

$$\begin{cases} -\partial_t u - \nu\Delta u + \frac{1}{2}|\nabla u|^2 & = \frac{1}{2}|x - \mu_m(t)|^2, \\ \ \partial_t m - \nu\Delta m - \operatorname{div}(m\nabla u) & = 0, \end{cases}$$

subject to the terminal and initial conditions:

$$u(\cdot, T) = 0, \quad m(\cdot, 0) = m_0^*.$$

The term $\mu_m(t)$ is defined as the first moment of the density function $m$, given by

$$[\mu_m(t)]_i = \int_{\mathbb{R}^d} y_i m(y,t)dy.$$

The initial density $m_0^*$ corresponds to a Gaussian distribution with mean $\mu_0 \in \mathbb{R}^d$ and a diagonal covariance matrix $\Sigma_0 \in \mathbb{R}^{d\times d}$.

In [7], the authors demonstrated that the unique classical solution to the HJB equation takes the quadratic form:

$$u(x,t) = \frac{1}{2}\Pi(t)|x|^2 + s(t)\cdot x + c(t),$$

where the coefficient functions are explicitly given by

$$\begin{aligned} \Pi(t) &= \tanh(T - t), \\ s(t) &= -\Pi(t)\mu_0, \\ c(t) &= \frac{1}{2}\Pi(t)|\mu_0|^2 + \nu\log(\cosh(T - t)). \end{aligned}$$

Furthermore, the density function remains Gaussian for all times $t$, with mean $\mu_0$ and a time-dependent diagonal covariance matrix $\Sigma(t)$ whose diagonal elements are

given by

$$(\Sigma(t))_{i,i} = (\Sigma_0)_{i,i}\frac{\cosh(T-t)^2}{\cosh(T)^2} + 2\nu\cosh(T-t)^2(\tanh(T) - \tanh(T-t)).$$

We note that these two classes of functions admit low-rank TT representations. In particular, for quadratic value functions it has been shown in [23, 40] that, to achieve an approximation error $\epsilon$, the required TT ranks grow only like $O(\log(\epsilon^{-1}))$. For weakly correlated Gaussian densities, [51] proves that the TT ranks increase at most poly-logarithmically with respect to the dimension $d$ and the inverse of the accuracy $1/\epsilon$.

We fix the dimension $d = 3$ and the space–time domain $Q = [-2.5, 2.5]^3 \times [0, T]$ with final time $T = 0.25$. We consider two examples: a deterministic case with $\nu = 0$ and a small-viscosity case with $\nu = 10^{-3}$. The initial mean is set to $\mu_0 = (0.1, 0.1, 0.1)$, and the diagonal elements of the initial covariance matrix $\Sigma_0$ are specified as $(\Sigma_0)_{i,i} = 0.5$ for all $i \in \{1, \dots, d\}$. For the numerical discretization, we employ $n_{\text{HJB}} = 3$ and $n_{\text{FP}} = 40$, which ensures a $L^2$ error of order $O(10^{-10})$ for the approximation of the initial density $m_0^*$. In this setting, we fix the stopping tolerance $\varepsilon_{\text{stop}} = 10^{-5}$ and the smoothing parameter to $\delta_n = 1$, i.e., no smoothing is applied, since it was not found to be necessary.

| $\Delta t$ | $E_2(u)$ | Order | $E_2(m)$ | Order | CPU (s) |
|---|---|---|---|---|---|
| $1.25 \times 10^{-1}$ | $3.71 \times 10^{-2}$ | – | $3.02 \times 10^{-2}$ | – | 1.66 |
| $6.25 \times 10^{-2}$ | $1.85 \times 10^{-2}$ | 1.01 | $1.49 \times 10^{-2}$ | 1.02 | 2.90 |
| $3.12 \times 10^{-2}$ | $9.10 \times 10^{-3}$ | 1.02 | $7.50 \times 10^{-3}$ | 0.99 | 5.95 |
| $1.56 \times 10^{-2}$ | $4.50 \times 10^{-3}$ | 1.02 | $3.70 \times 10^{-3}$ | 1.02 | 7.52 |

Table 11: Test 6.3, non-local MFG system. $L^2$ errors and CPU time for the SL1 scheme for the case $\nu = 0$, $d = 3$, with $n_{HJB} = 3$ and $n_{FP} = 40$.

| $\Delta t$ | $E_2(u)$ | Order | $E_2(m)$ | Order | CPU (s) |
|---|---|---|---|---|---|
| $1.25 \times 10^{-1}$ | $5.32 \times 10^{-3}$ | – | $2.23 \times 10^{-4}$ | – | 3.08 |
| $6.25 \times 10^{-2}$ | $1.31 \times 10^{-3}$ | 2.03 | $4.61 \times 10^{-5}$ | 2.27 | 6.24 |
| $3.12 \times 10^{-2}$ | $3.24 \times 10^{-4}$ | 2.01 | $1.06 \times 10^{-5}$ | 2.12 | 9.54 |
| $1.56 \times 10^{-2}$ | $8.07 \times 10^{-5}$ | 2.01 | $2.59 \times 10^{-6}$ | 2.04 | 17.82 |

Table 12: Test 6.3, non-local MFG system. $L^2$ errors and CPU time for the SL2 scheme for the case $\nu = 0$, $d = 3$, with $n_{HJB} = 3$ and $n_{FP} = 40$.

In Table 11 and 12 we show for several refinements of the time mesh, the CPU time, the errors $E_2(u)$ and $E_2(m)$ in the case $\nu = 0$, and the corresponding numerical convergence rate for the first order SL1 scheme and the second order SL scheme respectively. The first order SL1 scheme exhibits a first order convergence rate (Table 11), whereas the second order SL scheme achieves a second order convergence rate (Table 12), albeit at a higher per-step cost. It is worth noting that, for deterministic problems, the SL2p and SL2e schemes coincide; in this case, we refer to the scheme simply as SL2.

Introducing a small viscosity, $\nu = 10^{-3}$, Tables 13 and 14 show that the first order SL1 scheme maintains its first order convergence rate, while the second order

| $\Delta t$ | $E_2(u)$ | Order | $E_2(m)$ | Order | CPU (s) |
|---|---|---|---|---|---|
| $1.25 \times 10^{-1}$ | $3.71 \times 10^{-2}$ | – | $3.06 \times 10^{-2}$ | – | 5.15 |
| $6.25 \times 10^{-2}$ | $1.85 \times 10^{-2}$ | 1.00 | $1.54 \times 10^{-2}$ | 0.99 | 8.62 |
| $3.12 \times 10^{-2}$ | $9.10 \times 10^{-3}$ | 1.02 | $7.70 \times 10^{-3}$ | 1.00 | 13.86 |
| $1.56 \times 10^{-2}$ | $4.50 \times 10^{-3}$ | 1.02 | $3.90 \times 10^{-3}$ | 0.98 | 23.02 |

Table 13: Test 6.3 non-local MFG system. $L^2$ errors and CPU time for the SL1 scheme for the case $\nu = 10^{-3}$, $d = 3$, with $n_{HJB} = 3$ and $n_{FP} = 40$.

| $\Delta t$ | $E_2(u)$ | Order | $E_2(m)$ | Order | CPU (s) |
|---|---|---|---|---|---|
| $1.25 \times 10^{-1}$ | $5.16 \times 10^{-3}$ | – | $2.51 \times 10^{-4}$ | – | 9.91 |
| $6.25 \times 10^{-2}$ | $1.26 \times 10^{-3}$ | 2.04 | $6.12 \times 10^{-5}$ | 2.03 | 23.26 |
| $3.12 \times 10^{-2}$ | $3.08 \times 10^{-4}$ | 2.03 | $1.84 \times 10^{-5}$ | 1.74 | 35.91 |
| $1.56 \times 10^{-2}$ | $7.40 \times 10^{-5}$ | 2.06 | $6.48 \times 10^{-6}$ | 1.50 | 64.55 |

Table 14: Test 6.3 non-local MFG system. $L^2$ errors and CPU time for the SL2p scheme for the case $\nu = 10^{-3}$, $d = 3$, with $n_{HJB} = 3$ and $n_{FP} = 40$.

SL2p scheme achieves a near-quadratic rate in the HJB error. The FP error convergence rate, however, is mildly reduced at the finest resolutions. For comparable CPU times, however, the SL2p scheme attains errors smaller by two orders of magnitude, particularly for the FP equation, highlighting its superior efficiency and accuracy.

For all test cases, the TT rank associated with the value function remains equal to 3, while the maximal TT rank required for the Fokker–Planck solution is 14 for the SL2p scheme and 17 for the SL1 scheme, indicating the higher regularity of the high-order approximation.

It is also worth noting that the SL2e approximation produces errors of comparable magnitude while requiring additional computational effort; for this reason, its results are omitted.

**Preservation of Mass and First Moments.** An important feature of numerical methods for mean field games is their ability to preserve structural properties of the FP equation. In particular, we focus here on the preservation of the total mass of the density and its first moment. We denote by $M_T(m) = \int_{\mathbb{R}^d} m(x,T)\,dx$ the total mass of the exact solution of the FP equation at the final time $T$, and by $\widetilde{M}_T(m)$ the corresponding quantity computed with the TT approximant. Similarly, letting $\mu_m(t)$ denote the vector of first moments of the exact solution, we write $\widetilde{\mu}_m(t)$ for the moments obtained from the approximated solution. Based on these quantities, we define the errors $E_M = |M_T(m) - \widetilde{M}_T(m)|$ and $E_\mu = \|\mu_m(T) - \widetilde{\mu}_m(T)\|$.

Table 15 reports a comparison between the SL1 and SL2p schemes when the problem is computed with $\nu = 10^{-3}$. We observe that the second order approximation consistently delivers higher accuracy, with errors in both mass and first moments reaching magnitudes on the order of $10^{-7}$ to $10^{-6}$ for the finer discretizations. This confirms its advantage over the first order approach in preserving key integral quantities of the FP solution.

Overall, these results demonstrate that higher order schemes not only improve pointwise accuracy, but also enhance the preservation of physically relevant quantities such as mass and first moments.

| | SL1 scheme | | SL2p scheme | |
|---|---|---|---|---|
| $\Delta t$ | $E_M$ | $E_\mu$ | $E_M$ | $E_\mu$ |
| $1.25 \times 10^{-1}$ | $1.40 \times 10^{-3}$ | $1.30 \times 10^{-3}$ | $1.80 \times 10^{-4}$ | $1.75 \times 10^{-5}$ |
| $6.25 \times 10^{-2}$ | $1.20 \times 10^{-3}$ | $1.00 \times 10^{-3}$ | $2.23 \times 10^{-5}$ | $1.38 \times 10^{-6}$ |
| $3.12 \times 10^{-2}$ | $8.00 \times 10^{-4}$ | $6.00 \times 10^{-4}$ | $1.91 \times 10^{-6}$ | $7.00 \times 10^{-7}$ |
| $1.56 \times 10^{-2}$ | $4.00 \times 10^{-4}$ | $3.00 \times 10^{-4}$ | $8.23 \times 10^{-7}$ | $9.79 \times 10^{-7}$ |

Table 15: Test 6.3. Comparison of SL1 vs. SL2p schemes, in terms of errors in mass and first-momentum conservation in the case $\nu = 10^{-3}$, $d = 3$.

**High-dimensional deterministic case.** In this section, we analyze the behavior of the proposed methodology in the deterministic setting, i.e. when $\nu = 0$, varying the spatial dimension $d$. As already noted in Remark 2.1, in this case only one characteristic needs to be computed to evaluate the solution at a given point, since the dynamics follow the deterministic flow without any stochastic sampling. We investigate how this simplification is reflected in the scaling with respect to the dimension, by reporting the CPU times together with the corresponding exponential and polynomial fits.

We set $\Delta t = T/16$ for the SL1 scheme and $\Delta t = T/4$ for the SL2 scheme to achieve comparable $L^2$ errors in the HJB equation in any dimension $d$. Table 16 reports the relative $L^2$ errors and CPU times obtained with both methods as the dimension $d$ increases. The results show that the errors $E_2(u)$ and $E_2(m)$ remain nearly constant across all tested dimensions, while CPU times exhibit a moderate growth with $d$. We observe that the SL1 scheme is consistently slower than SL2. Although both methods yield similar accuracy in the HJB component, SL2 achieves errors in the FP equation that are approximately two orders of magnitude smaller. We observe that the TT rank associated with the value function remains fixed at 3, while the TT rank for the Fokker–Planck solution grows moderately from 12 in dimension $d = 3$ to 18 in dimension $d = 8$ for both methods. To further quantify this scaling behavior, Table 17 presents the fitted coefficients and corresponding $R^2$ values for both exponential and polynomial growth models. The analysis indicates that CPU times are best described by polynomial growth, with orders $O(d^{2.6})$ for SL2 and $O(d^{2.4})$ for SL1. This represents a notable improvement compared to the steeper $O(d^3)$ and $O(d^4)$ growth observed in the first numerical test, where diffusion was present.

| | SL1 | | | SL2 | | |
|---|---|---|---|---|---|---|
| $d$ | $E_2(u)$ | $E_2(m)$ | CPU (s) | $E_2(u)$ | $E_2(m)$ | CPU (s) |
| 3 | $2.10 \times 10^{-3}$ | $3.60 \times 10^{-3}$ | 6.73 | $1.31 \times 10^{-3}$ | $4.61 \times 10^{-5}$ | 3.33 |
| 4 | $2.10 \times 10^{-3}$ | $4.90 \times 10^{-3}$ | 19.17 | $1.31 \times 10^{-3}$ | $5.96 \times 10^{-5}$ | 7.66 |
| 5 | $2.10 \times 10^{-3}$ | $4.50 \times 10^{-3}$ | 31.17 | $1.31 \times 10^{-3}$ | $6.18 \times 10^{-5}$ | 11.66 |
| 6 | $2.10 \times 10^{-3}$ | $6.30 \times 10^{-3}$ | 43.83 | $1.31 \times 10^{-3}$ | $8.99 \times 10^{-5}$ | 21.33 |
| 7 | $2.10 \times 10^{-3}$ | $4.80 \times 10^{-3}$ | 56.32 | $1.31 \times 10^{-3}$ | $9.22 \times 10^{-5}$ | 31.03 |
| 8 | $2.10 \times 10^{-3}$ | $3.90 \times 10^{-3}$ | 83.46 | $1.31 \times 10^{-3}$ | $7.18 \times 10^{-5}$ | 45.12 |

Table 16: Test 6.3, non-local MFG system. Comparison of the SL1 and SL2 schemes in the case $\nu = 0$, varying the dimension $d$.

| Method | Fitting | Coefficient $a$ | Exponent $b$ | $R^2$ |
|---|---|---|---|---|
| SL1 | Exponential | 2.42 | 0.46 | 0.918 |
| SL1 | Polynomial | 0.56 | 2.43 | **0.981** |
| SL2p | Exponential | 0.87 | 0.51 | 0.961 |
| SL2p | Polynomial | 0.19 | 2.63 | **0.998** |

Table 17: Test 6.3, non-local MFG system. Results of the fittings for the SL1 and SL2 schemes using both exponential and power-law models for the case $\nu = 0$.

**6.4. 3D MFG with local linear coupling.** We conclude the numerical investigation with a three-dimensional Mean Field Game featuring a linear local coupling. The purpose of this test is to assess the performance of the proposed TT-based semi-Lagrangian method in a setting for which no explicit or closed-form solution is available. This experiment extends the previous numerical tests by departing from analytically tractable quadratic–Gaussian configurations. We consider the MFG system

$$\begin{cases} -\partial_t u - \nu\Delta u + \frac{|\nabla u|^2}{2} = \gamma\, m, \\ \partial_t m - \nu\Delta m - \nabla\cdot(m\nabla u) = 0, \end{cases}$$

on the domain $Q = [-1,1]^3 \times [0,T]$, with final time $T = 1$, periodic boundary conditions in space, and diffusion coefficient $\nu = 0.1$. In the following experiments, we vary the coupling strength $\gamma$. The terminal condition for the HJB equation is prescribed by the quadratic cost $u(x,T) = 0.1|x|^2$ which penalizes the distance from the origin. The initial density is taken as the superposition of two Gaussian bumps centered at $(-0.5,-0.5,-0.5)$ and $(0.5,0.5,0.5)$, each with variance 0.5. The objective of this test is to drive the population towards the origin while accounting for congestion effects, whose intensity is controlled by the parameter $\gamma$. In contrast to the Linear–Quadratic–Gaussian framework considered in previous experiments, the present configuration does not preserve a quadratic–Gaussian structure due to the presence of the coupling term $\gamma m$ on the right-hand side of the HJB equation.

Time discretization is performed using a uniform time step $\Delta t = 0.05$, and the SL2p scheme is employed. The value function and the density are approximated using $n_{HJB} = 3$ and $n_{FP} = 20$, respectively. The fixed-point iteration is terminated once the stopping criterion $\epsilon_{\text{stop}} = 10^{-3}$ is satisfied, and a smoothing parameter $\delta_n = 2/(2+n)$ is adopted, in accordance with the convergence results reported in [54].

Since no exact solution is available, accuracy is assessed through a consistency check based on the discrete HJB equation. More precisely, we evaluate the residual of the semi-Lagrangian scheme associated with the HJB equation at the initial time step. Let $u_k$ and $m_k$ denote the numerical approximations of the value function and the density at time $t_k$, respectively, obtained at the final policy iteration. For a given spatial point $x \in \mathbb{R}^3$, the pointwise residual of the HJB equation is then defined as

$$\mathcal{R}(x,u_k) := \left| u_k(x) - \sum_{\ell\in\mathcal{L}} \omega_\ell \Big( u_{k+1}(\Psi^+_{\ell,k}(x)) + \tfrac{\Delta t}{2} f_{k+1}(\Psi^+_{\ell,k}(x)) \Big) - \tfrac{\Delta t}{2} f_k(x) \right|,$$

where $\{\omega_\ell\}_{\ell\in\mathcal{L}}$ are the quadrature weights of the SL2p scheme, $\Psi^+_{\ell,k}(\cdot)$ denotes the forward stochastic characteristic, and $f_k(x)$ is the running cost evaluated at $m_k$, as defined in Section 5. This residual measures the defect in the discrete dynamic programming principle induced by the numerical approximation. In the numerical experiments, $\mathcal{R}(x, u_k)$ is evaluated on a coarse $20 \times 20 \times 20$ Cartesian grid at the first time step $k = 0$, for a total of $N_c = 20^3$ points. We report the mean squared residual and the uniform residual, namely

$$R_2(u_0) = \frac{1}{N_c} \sum_{k=1}^{N_c} \mathcal{R}(x_k, u_0)^2, \quad R_\infty(u_0) = \max_{1 \le k \le N_c} \mathcal{R}(x_k, u_0).$$

Small residual values indicate that the computed solution is an accurate fixed point of the fully discrete HJB–FP scheme. Table 2 reports the residuals for the cases $\gamma \in \{1, 2\}$. The residuals are notably small for $\gamma = 1$, whereas they increase for $\gamma = 2$, reflecting the increased difficulty of the problem due to the stronger congestion effect. Furthermore, we observe that the TT ranks of the value function remain equal to 4 at all time steps, while the maximal TT rank attained by the FP approximation is 25.

To further illustrate the qualitative behavior of the solution, we visualize the density at selected time instances. For the density, we employ three-dimensional isosurfaces corresponding to high quantiles of the mass distribution. Rather than plotting arbitrary level sets, these isosurfaces are constructed so as to enclose a prescribed fraction of the total mass (e.g., 96% and 98%), thereby providing a geometric characterization of the effective support of the population density.

In Figure 2, we report the density isosurfaces at different time steps for $\gamma = 1$ (top panels) and $\gamma = 2$ (bottom panels), together with the target position (red star) the center of the initial Gaussian density (black circles).

For $\gamma = 1$, the two populations merge and concentrate near the origin, whereas for larger coupling strength the density remains more dispersed: although it approaches the origin, the population tends to avoid excessive congestion, resulting in a more fragmented spatial distribution.

| $\gamma$ | $R_2(u_0)$ | $R_\infty(u_0)$ |
|---|---|---|
| 1 | $1.48 \times 10^{-5}$ | $4.27 \times 10^{-4}$ |
| 2 | $3.16 \times 10^{-5}$ | $1.50 \times 10^{-3}$ |

Table 18: Test 6.4, 3D MFG with linear coupling. Discrete residuals $R_2$ and $R_\infty$ for the SL2p scheme with $n_{HJB} = 3$ and $n_{FP} = 20$, for varying values of the coupling parameter $\gamma$.

**Conclusions.** This work presented a novel numerical framework for the approximation of high-dimensional Mean Field Games, built on the synergy between semi-Lagrangian time discretizations and Tensor-Train representations. The central objective was to mitigate the curse of dimensionality while retaining accuracy and stability in the discretization of the coupled HJB and FP equations.

A key contribution of this study lies in the design of a second order semi-Lagrangian scheme whose computational complexity grows quadratically with the dimension. This represents a significant advance over classical exponential-cost constructions, enabling accurate high-dimensional simulations without resorting to excessively large numbers of quadrature nodes. The proposed polynomial scheme achieves second order consistency in time while maintaining moderate growth in the number of characteristics, thus offering a practical compromise between efficiency and accuracy.

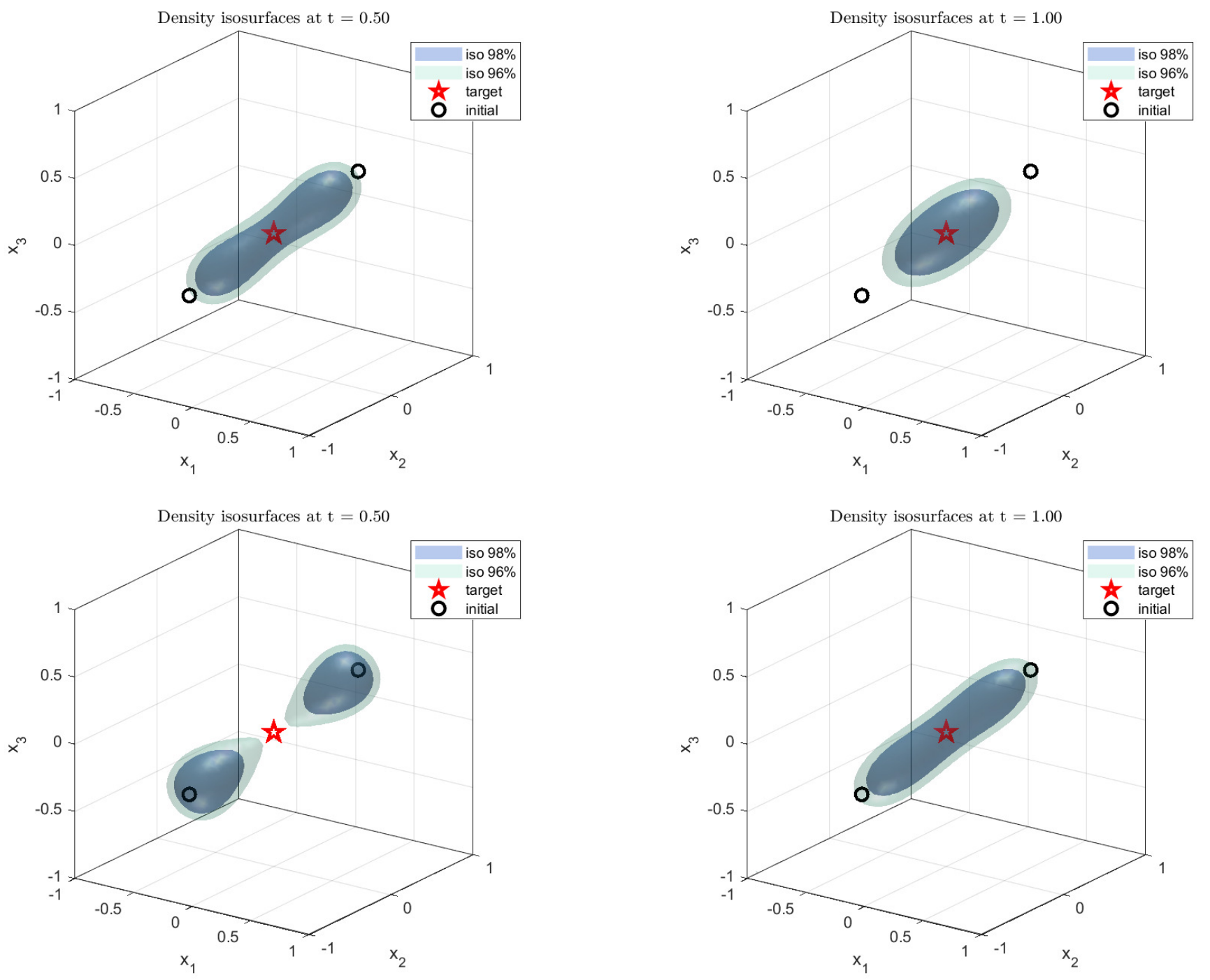

Fig. 2: Test 6.4, 3D MFG with linear coupling. FP density approximation at times $t = 0.5$ (left) and $t = 1$ (right), for $\gamma = 1$ (top) and $\gamma = 2$ (bottom).

Coupled with the Tensor-Train representation, which ensures compact storage and efficient evaluation of the solution and its derivatives, this strategy leads to a fully discrete algorithm that is both scalable and robust. Numerical experiments confirmed the theoretical convergence rates of the proposed schemes and highlighted their ability to preserve structural properties of the MFG system, such as asymptotic positivity and conservation of mass. In particular, the polynomial second order scheme was shown to be competitive with exponential schemes in low dimensions and preferable in higher dimensions, where its complexity scaling makes it markedly more efficient.

Future research directions include a rigorous theoretical analysis of the polynomial quadrature scheme, with attention to stability and error estimates when combined with TT truncations. Another direction is the extension of the framework to more complex applications, such as high-dimensional control problems arising in quadcopter dynamics and multi-agent robotic systems, where nonlinear dynamics and interaction effects play a significant role. These settings would allow the method to be assessed on problems that are closer to real-world scenarios.

## Appendix A. Proof of Proposition 3.1.

By expanding the solution $m$ in a Taylor series within the one-step error $\mathcal{T}_k(x)$, we prove the claim. Using the advection–diffusion equation (3.1), we obtain

$$\begin{aligned} m_t &= -q \cdot \nabla m + \nu \Delta m, \\ m_{tt} &= (q \cdot \nabla)^2 m - \nu q \cdot \nabla(\Delta m) - \nu \Delta(q \cdot \nabla m) + \nu^2 \Delta^2 m. \end{aligned}$$

Thus, fix $(x,t) \in \mathbb{R}^d \times [0, T-\Delta t]$. Observing that

$$\Delta(q \cdot \nabla m) = q \cdot \nabla(\Delta m) + \Delta q \cdot \nabla m + 2\sum_{i,j} \frac{\partial q_j}{\partial x_i} \frac{\partial^2 m}{\partial x_i \partial x_j},$$

the Taylor expansion of $m$ at $(x, t+\Delta t)$ reads
(A.1)

$$\begin{aligned} m(x,t+\Delta t) &= m + \Delta t\big(-q\cdot\nabla m + \nu \Delta m\big) + O(\Delta t^3) \\ &+ \frac{\Delta t^2}{2}\Big((q\cdot\nabla)^2 m - \nu(2q\cdot\nabla(\Delta m) + \Delta q \cdot \nabla m - 2\sum_{i,j}\frac{\partial q_j}{\partial x_i}\frac{\partial^2 m}{\partial x_i \partial x_j}) + \nu^2 \Delta^2 m\Big). \end{aligned}$$

In order to obtain the expansion of the right-hand side of (3.2), we proceed as follows. First, we rewrite the backward flow map (3.3) as

$$\Psi_\ell^-(x) = x + A_\ell(x).$$

Second, we observe that the displacement admits the following graded expansion:

$$A_\ell(x) = -\Delta t\, q(x) + \frac{\Delta t^2}{2} Dq(x)\, q(x) + \xi_\ell - \frac{\Delta t}{2} Dq(x)\, \xi_\ell - \frac{\Delta t}{4} \xi_\ell^\top D^2 q(x) \xi_\ell + R_\ell,$$

where $Dq$ and $D^2 q$ denote the Jacobian and Hessian of $q$, respectively, and $R_\ell = \mathcal{O}(\Delta t^3, \Delta t^2 |\xi_\ell|, \Delta t |\xi_\ell|^3)$. Finally, we perform a Taylor expansion of $m(x + A_\ell)$ up to fourth order with respect to the displacement $A_\ell$, which yields

$$\begin{aligned} m(x+A_\ell) = m(x) &+ \sum_{i=1}^d \frac{\partial m}{\partial x_i}(A_\ell)_i + \frac{1}{2}\sum_{i,j=1}^d \frac{\partial^2 m}{\partial x_i \partial x_j}(A_\ell)_i (A_\ell)_j \\ &+ \frac{1}{6}\sum_{i,j,k=1}^d \frac{\partial^3 m}{\partial x_i \partial x_j \partial x_k}(A_\ell)_i(A_\ell)_j(A_\ell)_k \\ &+ \frac{1}{24}\sum_{i,j,k,l=1}^d \frac{\partial^4 m}{\partial x_i \partial x_j \partial x_k \partial x_l}(A_\ell)_i(A_\ell)_j(A_\ell)_k(A_\ell)_l + \mathcal{O}(|A_\ell|^5). \end{aligned}$$

By the moment-matching conditions (3.5) on the quadrature nodes and weights, the weighted sums over $\ell$ reproduce the coefficients of the corresponding derivatives. The moment conditions, for $i, j, k, l \in \{1, \dots, d\}$, read

$$\sum_\ell w_\ell (A_\ell)_i = -\Delta t\, q_i + \frac{\Delta t^2}{2}\sum_j \left(q_j \frac{\partial q_i}{\partial x_j} - \nu \frac{\partial^2 q_j}{\partial x_i \partial x_j}\right) + \mathcal{O}(\Delta t^3),$$

$$\frac{1}{2}\sum_\ell w_\ell (A_\ell)_i (A_\ell)_j = \Delta t\, \nu\, \delta_{ij} + \frac{\Delta t^2}{2}\Big(q_i q_j - 2\nu \frac{\partial q_i}{\partial x_j}\Big) + \mathcal{O}(\Delta t^3),$$

$$\frac{1}{6}\sum_\ell w_\ell (A_\ell)_i (A_\ell)_j (A_\ell)_k = -\frac{\Delta t^2 \nu}{3}\big(q_i \delta_{jk} + q_j \delta_{ik} + q_k \delta_{ij}\big) + \mathcal{O}(\Delta t^3),$$

$$\frac{1}{24}\sum_\ell w_\ell (A_\ell)_i (A_\ell)_j (A_\ell)_k (A_\ell)_l = \frac{\Delta t^2 \nu^2}{6}\big(\delta_{ij}\delta_{kl} + \delta_{ik}\delta_{jl} + \delta_{il}\delta_{jk}\big) + \mathcal{O}(\Delta t^3).$$

Therefore, the scheme exactly reproduces the Taylor expansion of the solution up to order $\Delta t^2$. Hence, by substituting the time expansion of $m$ and the spatial expansion of the numerical update at $t = t_k$ into $\mathcal{T}_k(x)$, and exploiting the regularity of $m$ and $q$, it follows that $\mathcal{T}_k(x) = \mathcal{O}(\Delta t^3)$, uniformly in $x$, which proves the claim.

**Appendix B. Derivation of the Moment Matching System.**

In this section, we derive the coefficients $s$, $w_0$, $w_A$, and $w_D$ required for the construction of a second order approximation based on a polynomial quadrature rule. These coefficients are determined by enforcing the moment-matching conditions specified in equation (3.5). From the zeroth moment condition (mass conservation), summing the contributions of all nodes yields the equation:

$$w_0 + 2d\, w_A + 4\binom{d}{2} w_D = 1.$$

For the second moment condition, consider a fixed coordinate direction $x_i$. The axial nodes contribute through two points located at $\pm s e_i$, each with value $s^2$ and weight $w_A$. Additionally, for each pair $(i, j)$ with $j \neq i$, there are four face-diagonal nodes. In each of these nodes, the $i$th component of the vector $\xi$ is either $+s$ or $-s$, contributing $s^2$ to the second pure moment (due to squaring). For a fixed $i$, there are $d-1$ such pairs. This leads to the equation:

$$2w_A\, s^2 + 4(d-1)w_D\, s^2 = 2\nu\Delta t.$$

The condition for the pure fourth moment follows analogously. Since the contribution of each relevant node is $s^4$, we obtain:

$$2w_A\, s^4 + 4(d-1)w_D\, s^4 = 12(\nu\Delta t)^2.$$

For the mixed fourth moments, only the face-diagonal nodes contribute. For each pair $(i, j)$ with $i \neq j$, there are four such nodes, each contributing $s^4$. The resulting condition is:

$$4w_D\, s^4 = 4(\nu\Delta t)^2.$$

Solving this system of four equations yields the following expressions for the parameters:

$$s = \sqrt{6\nu\Delta t}, \quad \omega_0 = \frac{d^2 - 7d + 18}{18}, \quad \omega_A = \frac{4-d}{18}, \quad \omega_D = \frac{1}{36}.$$